\documentclass[a4paper]{article}
\usepackage{amsmath,amssymb,fancybox,longtable,graphics,amscd,multirow,theorem}
\usepackage[dvipdfm]{graphicx} 
\usepackage{caption,setspace}
\newtheorem{thm}{\textsc{Theorem}}[section]
\newtheorem{prop}{\textsc{Proposition}}[section]
\newtheorem{lem}{\textsc{Lemma}}[section]

\def\proof{\textsc{Proof. }}
\def\QED{$\Box$} 
\def\mbi#1{\boldsymbol{#1}} 

\def\Pic{\mathop{\mathrm{Pic}}\nolimits}

\def\rk{\mathop{\mathrm{rank}}\nolimits}
\def\sign{\mathop{\mathrm{sgn}}\nolimits}

\def\sing{\mathop{\mathrm{Sing}}\nolimits}

\theorembodyfont{\rmfamily}
\newtheorem{remark}{\textsc{Remark}}

\theorembodyfont{\rmfamily}

\begin{document} 
\title{Families of $K3$ surfaces and curves of $(2,3)$-torus type}
\author{Makiko Mase
\footnote{Tokyo Metropolitan University and University of Mannheim, mtmase@arion.ocn.ne.jp}
\footnote{{\it Keywords and phrases.} families of $K3$ surfaces that are double cover of the projective plane, curves of $(2,3)$-torus type, duality of Picard lattices, non-Galois triple covering of the projective plane }
}
\date{\empty}
\maketitle
\abstract{We study families of $K3$ surfaces obtained by double covering of the projective plane branching along curves of $(2,3)$-torus type. 
In the first part, we study the Picard lattices of the families, and a lattice duality of them. 
In the second part, we describe a deformation of singularities of Gorenstein $K3$ surfaces in these families.} \\

\noindent
{\it Running head. } Families of $K3$ surfaces and curves of $(2,3)$-torus type

\section{Introduction}
Plane sextic curves of $(2,3)$-torus type, which is defined by a polynomial of the form $F=F_2^3+F_3^2$ with polynomials $F_i$ of homogeneous degree $i$, that have at most simple singularities, which we simply call {\it curves of $(2,3)$-torus type}, are classified by Oka and Pho~\cite{OkaPho}, and Pho~\cite{Pho}. 

A compact complex connected algebraic surface $S$ is called a {\it Gorenstein $K3$ surface} if $S$ has at most simple singularities, has trivial canonical divisor, and the irregularity of $S$ is zero. 
If a Gorenstein $K3$ surface is non-singular, we simply call it a {\it $K3$ surface}. 
Note that for any Gorenstein $K3$ surface $S$, there exists a $K3$ surface $\tilde{S}$ such that $S$ and $\tilde{S}$ are birationally equivalent. 

Concerning algebraic curves on a $K3$ surface, the study of their Weierstrass semi-group and the existstance of curves that admit a given semi-group are quite authentic and intersting but there are not many known results: see for instance Komeda and Watanabe~\cite{KomedaWatanabe}. 
On the other hand, motivated originally in theoretical physics, many concepts have been  discovered  related to mirror symmetry for $2$-dimensional Calabi-Yau manifolds such as by Batyrev, Berglund and H$\ddot{\textnormal{u}}$bsch, Dolgachev, Ebeling and Takahashi, and Ebeling and Ploog~\cite{BatyrevMirror, BerglundHubsch, DolgachevMirror, EbelingTakahashi11, EbelingPloog}. 

In this article, we are to focus on families of $K3$ surfaces that are obtained as double covering of the projective plane branching along curves of $(2,3)$-torus type, and their Picard lattices.

Consider the generic non-Galois triple covering $X$ of $\mathbb{P}^2$ branching along a curve of $(2,3)$-torus type $B$, and then take the Galois closure $\hat{X}$ of $X$. 
Moreover, if $Z$ is the minimal model of the double covering $D(X\slash \mathbb{P}^2)$ of $\mathbb{P}^2$ branching along $B$, then, $\hat{X}$ is also obtained as the cyclic triple covering of $Z$. 

By this construction, and only in this case, one obtains a cyclic triple covering $\hat{X}\to Z$, where $Z$ is a $K3$ surface, and $\hat{X}$ is either an abelian surface (possibly with singularities), or a Gorenstein $K3$ surface. 
These two cases are distinguished by an invariant $\delta$ introduced by Ishida and Tokunaga~\cite{IshidaTokunaga}, in which, the latter case where $\delta=6$ is studied. 
The former case, in which one has $\delta=9$ is in fact studied in detail by Barth~\cite{Barth98} as an analogy of Nikulin's result~\cite{Nikulin75}. 
It is easily seen by the classification~\cite{OkaPho, Pho}, that $\delta=9$ occurs if and only if $\sing{B}=9A_2$; and together with the fact that the order of the fundamental group of the singularities should be divided by $3$, that $\delta=6$ occurs when $\sing{B}$ is one of the followngs: 
\[
\begin{array}{llll}
A_{17}, \, 2A_8, \,  3A_5, & A_2+A_{14}, & A_5 + A_{11}, & A_{11}+E_6, \\
3E_6,\, 6A_2, & 2A_2+A_{11}, & A_2+A_5+A_8, & 3A_2+A_8, \\
2A_2+2A_5, & 4A_2+A_5, & 2A_5+E_6, & A_5+2E_6, \\
2A_2+A_5+E_6, & 2A_2+2E_6,  & A_2+A_{8}+E_6, & 4A_2+E_6. 
\end{array}
\] 

\begin{remark}
The above list of $\sing{B}$ picked up in~\cite{Pho} covers all possible cases: indeed, with the aid of Proposition 1.1~\cite{IshidaTokunaga}. 

\end{remark}

The dual curve of a plane smooth cubic curve is a typical example (see \cite{BrieKnor}) of a curve of $(2,3)$-torus type $B$ with nine singularities of type $A_2$ (cusps). 


Since $\delta=9$ case is well-understood, we focus on the other case in this article. 

The defining equation of the Gorenstein $K3$ surface $D(X\slash \mathbb{P}^2)$ is given by 
\[
W^2 - F(X,Y,Z) = 0,
\]
where $F$ is the defining polynomial of the branch curve $B$. 
Being parameterised by the coefficients of the monomials in $W^2 - F(X,Y,Z)$, one can construct families of (Gorenstein) $K3$ surfaces. 
Such a family should be a subfamily of $K3$ surfaces parametrised by the complete anticanonical linear system of the weighted projective space $\mathbb{P}(1,1,1,3)$ with weights $(1,1,1,3)$, which is one of $95$ weights corresponding to simple $K3$ singularities classified by Yonemura~\cite{Yonemura}. 

The aim of this article is to study these families: 
we first construct polytopes $\Delta_1,\, \Delta_2$, and $\Delta_3$ such that the complete anticanonical linear systems of toric Fano $3$-folds associated to them parametrise families $\mathcal{F}_1,\,\mathcal{F}_2,\,\mathcal{F}_3$ of $K3$ surfaces obtained by a double covering of $\mathbb{P}^2$ branching along curves of $(2,3)$-torus type. 
Then, we shall give the Picard lattice of the families using toric geometry. 
More precisely, our main theorem is stated as follows: 

\noindent
{\bf Theorem \ref{MainThm}}\quad 
{\it 
{\rm (1)}\, The Picard lattices of the families $\mathcal{F}_1$ and $\mathcal{F}_2$ are respectively isometric to  $U\oplus \langle -2\rangle \oplus \langle -4\rangle$, and $U\oplus A_5$. \\
{\rm (2)}\, The Picard lattice of the family $\mathcal{F}_3$ is isometric to  $\langle -2\rangle \oplus \langle 2\rangle$. \\
{\rm (3)}\, The Picard lattice $\Pic_{\!\Delta_3}$ of the family $\mathcal{F}_3$ satisfies the duality 
\[
U\oplus \Pic_{\!\Delta_3^*} \simeq (\Pic_{\!\Delta_3})^\perp_{\Lambda_{K3}},
\]
where $\Delta_3^*$ is the polar dual of $\Delta_3$, $U$ is the  hyperbolic lattice of rank $2$, and $(L)^\perp_{\Lambda_{K3}} $ is the  orthogonal complement of a primitive sublattice $L$ in the $K3$ lattice $\Lambda_{K3}$.  } \\

We review fundamental facts of toric geometry necessary in this article in $\S 2$. 
The main theorem is proced in $\S 3$ after verifying invariants. 
In $\S 4$, we describe families that contain $K3$ surfaces obtained as double covering of $\mathbb{P}^2$ branching along curves of $(2,3)$-torus type. \\

For a curve $B$ of $(2,3)$-torus type, when $\sing(B)=A$ is the set of singularities of $B$, we call $[A]$ for the curve $B$. 
If $\sing(B)$ contains more than one singularities, we denote by $A=A' + A'', \, A=2A'=A'+A'$ {\it etc}. 

The singularities of type $A_n$ is the singularity given locally by $x^2+y^{n+1}=0$ for $n\geq 1$, and of type $E_6$ is given locally by $x^3+y^4=0$. 
In un-necessarily confusing way, we also mean by $A_n$ the root lattice of type $A_n$. 
The hyperbolic lattice of rank $2$ is denoted by $U$, and negative definite root 
lattice of rank $8$ is denoted by $E_8$. 
The {\it $K3$ lattice} is defined by $\Lambda_{K3} := U^{\oplus 3}\oplus E_8^{\oplus 2}$. 

\begin{ackn}
\textnormal{
The author thanks to Professor J.Komeda who gave her an opportunity to study Weierstrass semi-groups, and hopes that the result may enrich the study of that of algebraic curves on  $K3$ surfaces. 
Thanks to Professor C.Hertling for reading through the first draft. }
\end{ackn}

\section{Setups}
A {\it lattice} is a finitely-generated $\mathbb{Z}$-module with a non-degenerate bilinear form. 
Let $M$ be a lattice of rank $3$, and $N:={\rm Hom}_{\mathbb{Z}}(M,\mathbb{Z})$ be its dual. 

See Section 3.3 of~\cite{FultonToric} for basic facts on toric geometry. 

A {\it polytope} is a convex hull of finitely-many points, say, $v_1,\ldots, v_r$, in $M\otimes_{\mathbb{Z}}\mathbb{R}$, and is denoted by
\[
\Delta = {\rm Conv}\{v_1,\ldots, v_r\}. 
\]
For a polytope $\Delta$, a {\it vertex} is a $0$-dimensional face, an {\it edge} is a $1$-dimensional face, and a {\it face} is a $2$-dimensional face.  
A polytope is called {\it integral} if every vertex is in $M$. 

Denote by $\langle \cdot , \cdot \rangle : N\times M \to \mathbb{Z}$ a natural paring, which is the inner product in $\mathbb{R}^3$. 
Let $\Delta$ be a polytope in $M_{\mathbb{R}}$. 
Define the {\it polar dual} polytope of $\Delta$ by 
\[
\Delta^* := \left\{ y\in N\otimes_{\mathbb{Z}}\mathbb{R}\, |\, \langle y, x \rangle\geq -1 \, \textnormal{ for all }  x\in \Delta\right\}. 
\]

An integral polytope $\Delta$ is {\it reflexive} if $\Delta$ contains the origin in the interior as its only lattice point, and the polar dual $\Delta^*$ is also integral. 

It is known~\cite{BatyrevMirror} that a polytope $\Delta$ is reflexive if and only if the corresponding toric variety $\mathbb{P}_{\!\Delta}$ is Fano, in particular, the complete anticanonical linear system of $\mathbb{P}_{\!\Delta}$ parametrises a family $\mathcal{F}_{\!\Delta}$ of which general sections are Gorenstein $K3$ surfaces. 
Since every Gorenstein $K3$ surface is birationally equivalent to a unique $K3$ surface, we identify these two surfaces. 

Suppose $\Delta$ is a reflexive polytope, and take a general anticanonical section $Z$ of $\mathbb{P}_{\!\Delta}$.

It is also known~\cite{BatyrevMirror} that there exists a MPCP desingularization, which is achieved by a simultaneous toric desingularization, and in consequence, we obtain the resulting desingularized varieties $\widetilde{\mathbb{P}_{\!\Delta}}$ and $\widetilde{Z}$ of the ambient space $\mathbb{P}_{\!\Delta}$ and of $Z$, respectively. 

Denote by $\widetilde{\mathbb{P}_{\!\Delta}}$ and $\widetilde{Z}$ the resulting desingularized varieties. 
Define a restriction map $r:H^{1,1}(\widetilde{\mathbb{P}_{\!\Delta}}) \to H^{1,1}(\widetilde{Z})$ of Hodge $(1,1)$-parts. 
Note that the map $r$ is not necessarily surjective. 
Define 
\[
L_D(\Delta) :={\rm Im}\,{r}\cap H^2(\tilde{Z},\,\mathbb{Z}), \qquad L_0(\Delta):=(L_D(\Delta))^\perp_{H^2(\tilde{Z},\,\mathbb{Z})}. 
\]
Then, 
\[
L_0(\Delta)={\rm coker}\,{r}\cap H^2(\tilde{Z},\,\mathbb{Z}). 
\]
We call the rank of $L_0(\Delta)$ the {\it toric correction term}. 

Denote by $v_i,\, i=1,\ldots ,r$ the vectors starting from the origin and ending at vertices $e_i$ of $\Delta$. 
The vectors $v_i$ being as one-simplices, one can construct a fan $\Sigma=\Sigma(\Delta)$ associated to the polytope $\Delta$. 
It is easily observed that $\rk L_0(\Delta)=0$ if and only if $\Delta$ is {\it simplicial}, namely, any three of one-simplices of $\Sigma$ generate $M$, and it is also equivalent that the toric variety $\mathbb{P}_{\!\Delta}$ has at most orbifold singularities. 

Define the {\it Picard lattice of the family } $\mathcal{F}_{\!\Delta}$ to be the Picard lattice of $K3$ surfaces that are minimal models of any generic sections in $\mathcal{F}_{\!\Delta}$, which is known to be well-defined~\cite{Bruzzo-Grassi}. 
In other words, the Picard lattice of the family is a lattice generated by the irreducible components of the restrictions to ${-}K_{\!\mathbb{P}_{\!\Delta}}$ of generators, which are torus-action invarant, of the Picard group of the toric $3$-fold $\mathbb{P}_{\!\Delta}$. 
We denote it by $\Pic_{\!\Delta}$ and $\rho_{\!\Delta}$ be its rank. 

The toric correction term $\rk L_0(\Delta)$, $\rho_{\!\Delta}$, and intersection numbers $D.D'$ for restricted torus-invariant divisors $D, D'$ in $\Pic_{\!\Delta}$ can be combinatorically computed~\cite{Kobayashi}. 

For a face $\Gamma$ of any dimensional of $\Delta$, denote by $l(\Gamma)$ the number of lattice points in  $\Gamma$, and $l^*(\Gamma)$ the number of inner lattice points in  $\Gamma$. 

For an edge $\Gamma$ of $\Delta$, denote by $\Gamma^*$ the dual edge of $\Gamma$ in $\Delta^*$. 
The toric correction term is computed by
\begin{equation}
\rk L_0(\Delta) = \sum_{\Gamma: \textnormal{edge of } \Delta} l^*(\Gamma)l^*(\Gamma^*). \label{L0}
\end{equation}

The rank $\rho_{\!\Delta}$ is computed by
\begin{equation}
\rho_{\!\Delta} =  \sum_{\Gamma : \textnormal{edges of } \Delta}l(\Gamma^*) + \rk L_0(\Delta)-3. \label{Rho}
\end{equation}

A primitive vector $v$ that generates a ray of the fan defining the $3$-fold $\mathbb{P}_{\!\Delta}$ determines a torus-invariant divisor $\overline{{\rm Orb}(\mathbb{R}_{\geq 0}v)}$ on $\mathbb{P}_{\!\Delta}$. 
By construction, it is equivalent to take a lattice point $v^*$ in $\Delta^*$. 
The dual of $v^*$ is a face $F$ in $\Delta$. 
It is well-known that the number of lattice points in the interior of $F$ is equal to the genus of a smooth curve corresponding to the divisor. 
For $D:=\overline{{\rm Orb}(\mathbb{R}_{\geq 0}v)}|_{{-}K_{\!\mathbb{P}_{\!\Delta}}}$, one has
\begin{equation}
D.D=D^2 = 2l^*(F) -2. \label{SelfIntersection}
\end{equation}

Let $D$ and $D'$ be the restriction to ${-}K_{\!\mathbb{P}_{\!\Delta}}$ of torus-invariant divisors on $\mathbb{P}_{\!\Delta}$ corresponding to vertices $v$ and $v'$, and let $\Gamma$ be the edge of $\Delta$ which connects the vertices $v$ and $v'$. 
The intersection number is thus obtained by
\begin{equation}
D.D' = l^*(\Gamma)+1. \label{Intersection}
\end{equation}

\section{Main Results}
Define a lattice $M$ by 
\[
M:=\left\{ (a_0,\, a_1,\, a_2, \, a_3)\in\mathbb{Z}^4 \, | \, a_0+a_1+a_2+3a_3 \equiv 0 (6) \right\}. 
\]
The lattice $M$ is of  rank $3$, and one can take a basis $\{ e_1, e_2, e_3\}$ for $M$ over $\mathbb{Z}$, where 
\[
e_1 = (1, 0, -1, 0), \quad 
e_2 = (0, 1, -1, 0), \quad 
e_3 = (0, 0, -3, 1). 
\]
We can associate $M$ with the set of monomials $\mathcal{M}_6$ of weighted degree $6$ with weights $(1,1,1,3)$ by 
\[
\begin{array}{ccc}
M & \rightarrow & \mathcal{M}_6 \\
(a_0,\, a_1,\, a_2,\, a_3) & \mapsto & X^{a_0+1}Y^{a_1+1}Z^{a_2+1}W^{a_3+1}, 
\end{array}
\]
where the weights of $X, Y, Z, W$ are respectively $1,1,1,3$, and thus there exists a correspondence between lattice points in $M$ and such monomials. 

One can embed a polytope of which lattice points are labelled by monomials in $\mathcal{M}_6$ into $\mathbb{R}^3$ : 
express elements in $M$ as a linear combination of the chosen basis $\{ e_1,\, e_2,\, e_3\}$ of which the coefficients form a point in $\mathbb{R}^3$. 
Thus, one gets a correspondence between monomials and points in $\mathbb{R}^3$ under this choice of a basis, for which some examples are presented in Table \ref{MonomialLatticePoints}. 
\[
\begin{array}{ccc}
W^2 \leftrightarrow (-1, -1, 1), & Y^6 \leftrightarrow (-1, 5, -1), & X^6 \leftrightarrow (5, -1, -1), \\
Z^6\leftrightarrow (-1,-1,-1), & Y^3Z^3 \leftrightarrow (-1, 2, -1), & Y^2Z^4 \leftrightarrow (-1, 1, -1), \\
X^4Z^2 \leftrightarrow (3, -1, -1), & X^2Z^4 \leftrightarrow (1, -1, -1). 
\end{array}
\]
\begingroup
\captionof{table}{Correspondence of monomials and points in $M$. }\label{MonomialLatticePoints}
\endgroup

\bigskip

Define integral polytopes $\Delta_1,\, \Delta_2$, and $\Delta_3$ by 
\[
\Delta_1 := {\rm Conv}\left\{ (-1, -1, 1),\, (-1, 1, -1),\, (3, -1, -1),\, (5, -1, -1),\, (-1, 5, -1)\right\}, 
\]
\[
\Delta_2 := {\rm Conv}\left\{ (-1, -1, 1),\, (-1, 2, -1),\, (3, -1, -1),\, (5, -1, -1),\, (-1, 5, -1)\right\}, 
\]
and 
\[
\Delta_3 := {\rm Conv}\left\{  (-1, -1, 1),\, (-1, 1, -1),\, (1, -1, -1),\, (5, -1, -1),\, (-1, 5, -1)\right\} ,
\]
respectively. 
\begin{figure}[htb!] 
\begin{center}
\includegraphics[width=1.0\linewidth]{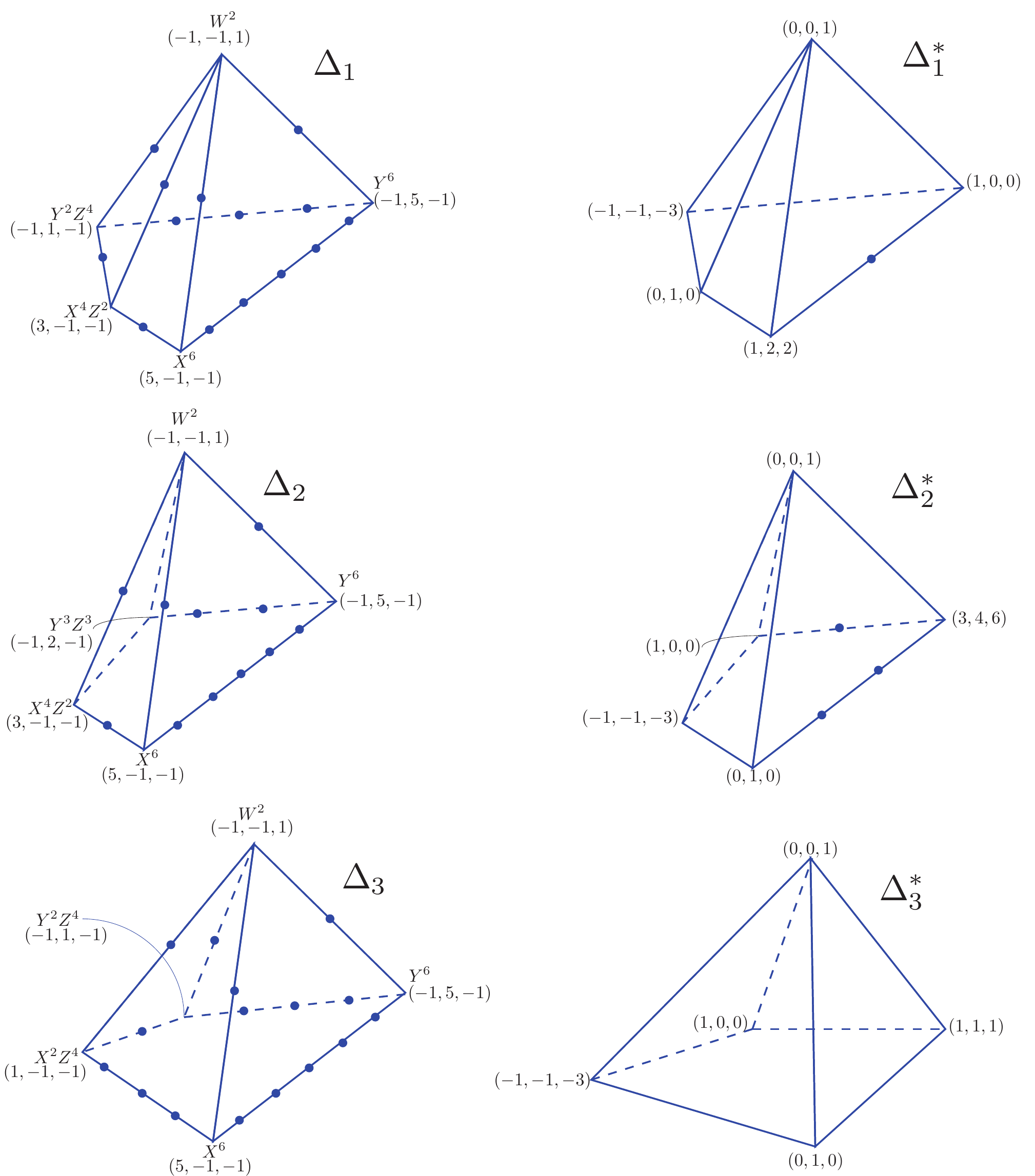} 
\end{center}
\caption{Polytopes $\Delta_1,\,\Delta_2,\,\Delta_3,\,\Delta_1^*,\,\Delta_2^*$ and $\Delta_3^*$. }
\label{Polytopes123}
\end{figure}%
See Figure~\ref{Polytopes123}. 

\begin{prop}
The polytopes $\Delta_1,\, \Delta_2$, and $\Delta_3$ are reflexive. 
\end{prop}
\proof
It is clear that the origin is the only lattice point contained in the polytopes. 
By a direct computation, the polar dual polytopes $\Delta_1^*,\, \Delta_2^*$, and $\Delta_3^*$ of $\Delta_1,\, \Delta_2$, and $\Delta_3$ are as follows: 
\[
\Delta_1^* = {\rm Conv}\left\{ (0, 0, 1),\, (-1, -1, -3),\, (0, 1, 0),\, (1, 2, 2),\, (1, 0, 0)\right\},
\]
\[
\Delta_2^* = {\rm Conv}\left\{ (0, 0, 1),\, (-1, -1, -3),\, (0, 1, 0),\, (3, 4, 6),\, (1, 0, 0)\right\},
\]
and 
\[
\Delta_3^* = {\rm Conv}\left\{ (0, 0, 1),\, (-1, -1, -3),\, (0, 1, 0),\, (1, 1, 1),\, (1, 0, 0)\right\},
\]
resepctively. 
Thus $\Delta_1^*,\, \Delta_2^*$, $\Delta_3^*$ are integral as well. 
\QED

\begin{prop}\label{L0Rho}
We have $\rk L_0(\Delta_1)=1$,\, $\rk L_0(\Delta_2)=2$, and $\rk L_0(\Delta_3)=0$, and $\rho_{\!\Delta_1}=4$,\, $\rho_{\!\Delta_2} = 7$,\, $\rho_{\!\Delta_3}=2$, $\rho_{\!\Delta_1^*}=17$,\, $\rho_{\!\Delta_2^*} = 15$, and $\rho_{\!\Delta_3^*}=18$. 
In particular, $\rho_{\!\Delta_i}+\rho_{\!\Delta_i^*}>20$ for $i=1,2$, and $\rho_{\!\Delta_3}+\rho_{\!\Delta_3^*}=20$ hold. 
\end{prop}
\proof
There exists a lattice point on the edge
\[
\Gamma_1 = {\rm Conv}\{ (-1, -1, 1),\, (-1, 1, -1)\}
\]
of $\Delta_1$, and a lattice point on its dual edge
\[
\Gamma_1^* = {\rm Conv}\{ (1, 0, 0),\, (1, 2, 2)\}
\]
of $\Delta_1^*$. 
There is no more edge on $\Delta_1$ that contributes $\rk L_0(\Delta_1)$. 
Thus, by the formula (\ref{L0}), 
\[
\rk L_0(\Delta_1) = l^*(\Gamma_1)l^*(\Gamma_1^*) = 1\cdot 1 = 1. 
\] 
By the formula (\ref{Rho}), one has $\rho_{\!\Delta_1}=(5+1)+1-3 =4$, and $\rho_{\!\Delta_1^*}=19+1-3 =17$. 
Clearly, $\rho_{\!\Delta_1}+\rho_{\!\Delta_1^*}=4+17=21>20$. 

There exists a lattice point on the edge 
\[
\Gamma_2 = {\rm Conv}\{ (-1, -1, 1),\, (3, -1, -1)\}
\]
of $\Delta_2$, and two lattice points on its dual edge
\[
\Gamma_2^* = {\rm Conv}\{ (3, 4, 6),\, (0, 1, 0)\}
\]
of $\Delta_2^*$. 
There is no more edge on $\Delta_2$ that contributes $\rk L_0(\Delta_2)$. 
Thus, by the formula (\ref{L0}), 
\[
\rk L_0(\Delta_2) = l^*(\Gamma_2)l^*(\Gamma_2^*) = 1\cdot 2 = 2. 
\] 
By the formula (\ref{Rho}), one has $\rho_{\!\Delta_2}=5+3+2-3 =7$, and $\rho_{\!\Delta_2^*}=16+2-3 =15$. 
Clearly, $\rho_{\!\Delta_2}+\rho_{\!\Delta_2^*}=7+15=22>20$. 

There does not exist an edge on $\Delta_3$ that contributes $\rk L_0(\Delta_3)$. 
Thus, by the formula (\ref{L0}), $\rk L_0(\Delta_3) = 0$. 
By the formula (\ref{Rho}), one has $\rho_{\!\Delta_3}=5-3 =2$, and $\rho_{\!\Delta_3^*}=21-3 =18$. 
Clearly, $\rho_{\!\Delta_3}+\rho_{\!\Delta_3^*}=2+18=20$. 
\QED

\bigskip

Denote by $\mathcal{F}_1,\, \mathcal{F}_2,\, \mathcal{F}_3$, and $\mathcal{F}_3'$ the families of $K3$ surfaces parametrised by the complete anticanonical linear systems of toric Fano $3$-folds $\mathbb{P}_{\!\Delta_1}$,\, $\mathbb{P}_{\!\Delta_2}$,\, $\mathbb{P}_{\!\Delta_3}$, and $\mathbb{P}_{\!\Delta_3^*}$, respectively. 
Here, $\Delta_3^*$ is the polar dual polytope of $\Delta_3$. 
Denote by $\Pic_{\Delta}$ the Picard lattice of the family of $K3$ surfaces that is associated to a reflexive polytope $\Delta$. 

\begin{remark}
We have seen that the sum of Picard numbers $\rho_{\!\Delta_3}$ and $\rho_{\!\Delta_3^*}$ coincides with the rank of the unimodular lattice $U^{\oplus 2}\oplus E_8^{\oplus 2}$, and that the rank of $L_0(\Delta_3)$ is $0$ means that the toric $3$-fold $\mathbb{P}_{\!\Delta_3}$ is simplicial. 
In~\cite{Mase17}, it is concluded that if a toric Fano $3$-fold $\mathbb{P}_{\!\Delta}$ is simplicial, then, the family $\mathcal{F}_{\Delta}$ of $K3$ surfaces is lattice dual in the sense that 
\[
\left(\Pic_{\!\Delta}\right)_{\Lambda_{K3}}^\perp \simeq U\oplus\Pic_{\!\Delta^*}
\]
holds. 
Thus, in our situation here, we expect that the family $\mathcal{F}_3$ constructed by the toric $3$-fold $\mathbb{P}_{\!\Delta}$ might be lattice dual. 
Therefore, we have to study the dual $\Delta_3^*$, or equivalently, the family $\mathcal{F}_3'$ constructed by the toric $3$-fold $\mathbb{P}_{\!\Delta^*}$. 
\end{remark}

\begin{remark}
Let $\Delta$ be a reflexive polytope. 
We occasionally identify the complete anticanonical linear system of the toric Fano $3$-fold $\mathbb{P}_{\!\Delta}$, and the family $\mathcal{F}_{\Delta}$. 
Indeed, a section $f\in |{-}K_{\mathbb{P}_{\!\Delta}}|$ determines a surface $(f=0)$ in $\mathcal{F}_{\Delta}$. 
Thus, we may also call $(f=0)$ a section as long as there is no confusion. 
\end{remark}

\begin{figure}[!htb]
\begin{center}
\includegraphics[width=.4\linewidth]{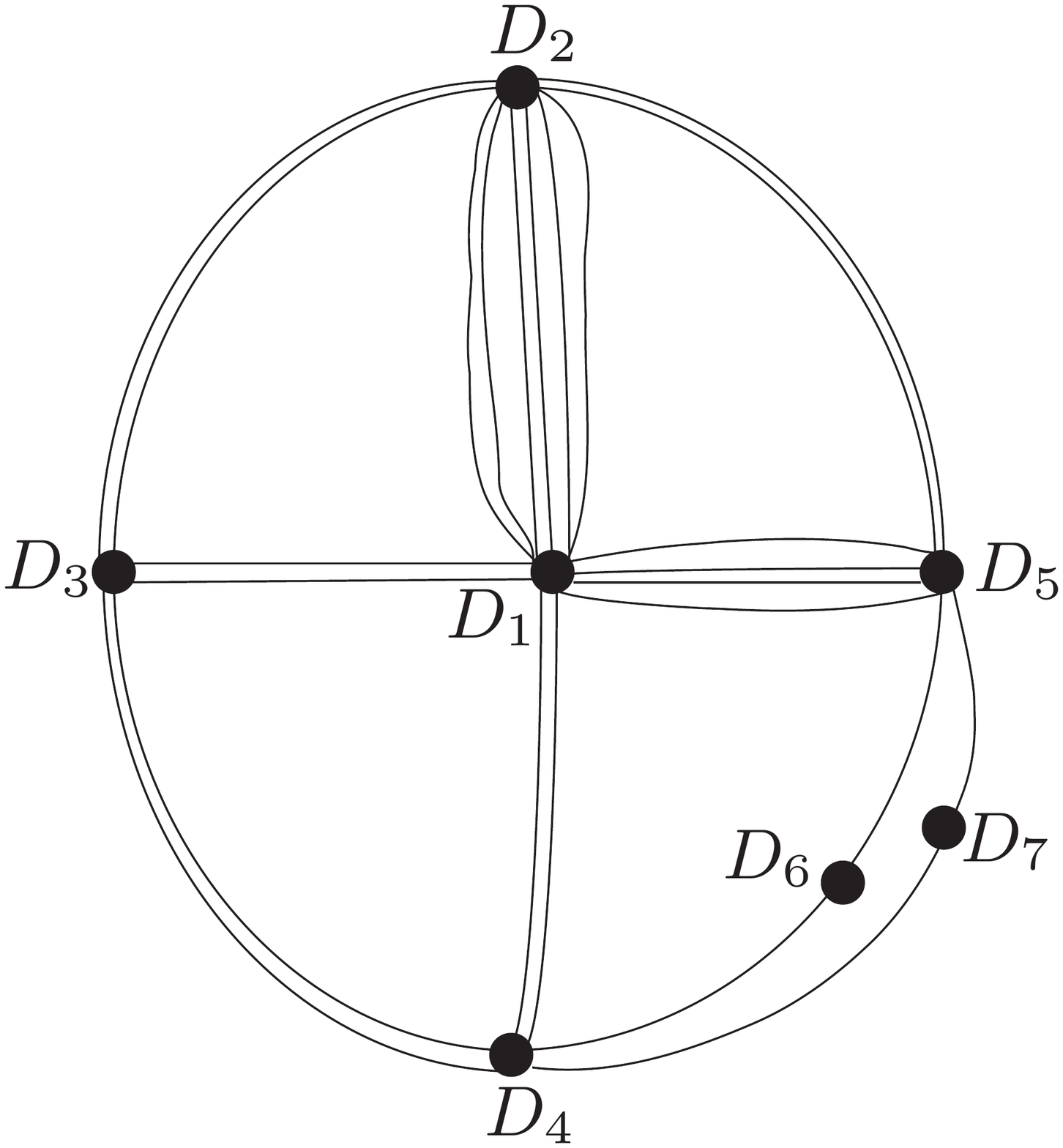}
\caption{Nodes are divisors on a general section in $\mathcal{F}_1$. 
There are $n$ lines connecting $D_i$ and $D_j$ if $D_i.D_j=n$. }
\label{Case1}
\end{center}
\end{figure}

\begin{thm}\label{MainThm}
{\rm (1)}\, The Picard lattices of the families $\mathcal{F}_1$ and $\mathcal{F}_2$ are respectively isometric to  $U\oplus \langle -2\rangle \oplus \langle -4\rangle$, and $U\oplus A_5$. \\
{\rm (2)}\, The Picard lattice of the family $\mathcal{F}_3$ is isometric to  $\langle -2\rangle \oplus \langle 2\rangle$. \\
{\rm (3)}\, The Picard lattice $\Pic_{\!\Delta_3}$ of the family $\mathcal{F}_3$ satisfies the duality 
\[
U\oplus \Pic_{\!\Delta_3^*} \simeq (\Pic_{\!\Delta_3})^\perp_{\Lambda_{K3}},
\]
where $\Delta_3^*$ is the polar dual of $\Delta_3$, $U$ is the  hyperbolic lattice of rank $2$, and $(L)^\perp_{\Lambda_{K3}} $ is the  orthogonal complement of a primitive sublattice $L$ in the $K3$ lattice $\Lambda_{K3}$. 
\end{thm}
\proof
(1)\, We lebel the primitive vectors that generate rays of the fan defining the $3$-fold $\mathbb{P}_{\!\Delta_1}$, or equivalently, the lattice points in $\Delta_1^*$ as follows:
\[
\begin{array}{lll}
v_1=(0,0,1), & v_2=(-1,-1,-3), & v_3=(0,1,0),\\
v_4=(1,2,2), & v_5=(1,0,0), & v_6=(1,1,1). 
\end{array}
\]
Let $D_i:=\overline{{\rm Orb}(\mathbb{R}_{\geq 0}v_i)}|_{{-}K_{\mathbb{P}_{\!\Delta_1}}}$ for $i=1,\ldots ,5$ be restricted torus-invatiant divisors, and $D_6, \, D_7$ be components of the divisor $\overline{{\rm Orb}(\mathbb{R}_{\geq 0}v_6)}|_{{-}K_{\mathbb{P}_{\!\Delta_1}}}$. 
One computes the self-intersection numbers by the formula (\ref{SelfIntersection}), 
\[
\begin{array}{lll}
D_1^2 = 2\cdot 8 -2 = 14, &
D_2^2 = 2\cdot 2 -2 = 2, &
D_3^2 = 2\cdot 0 -2 = -2, \\
D_4^2 = 2\cdot 0 -2 = -2, &
D_5^2 = 2\cdot 1 -2 = 0, &
D_6^2 = D_7^2 = -2. 
\end{array}
\]
One also has the graph of intersections among these divisors by the formula (\ref{Intersection}) as in Figure \ref{Case1}.  

By solving the linear system
\[
\sum_{i=1}^7 ({}^tv_i, \mbi{e}_j)D_i = 0, 
\]
for $j=1,2,3$, where $\mbi{e}_j$ is the $j$-th column of the identity matrix of size $3$, and $(\mbi{x}, \mbi{y})$ is the inner product on $\mathbb{R}^3$,  one obtains a set of linearly-independent divisors $\mathcal{B}=\{ D_4,\, D_5,\, D_6,\, D_7\}$, the intersection matrix with respect to which is 
\[
A_{\mathcal{B}} = \begin{pmatrix} -2 & 0 & 1 & 1\\ 0 & 0 & 1 & 1 \\ 1 & 1 & -2 & 0 \\ 1 & 1 & 0 & -2\end{pmatrix}. 
\]
By a translation $PA_{\mathcal{B}}{}^tP$ with 
\[
P = \begin{pmatrix} 0 & 1 & 1 & 0 \\ 0 & 1 & 0 & 0 \\ 1 & -1 & 0 & 0 \\ 0 & -2 & -1 & 1 \end{pmatrix}, 
\]
one gets a new basis 
\[
\mathcal{B}' = \{ D_5+D_6,\, D_5,\, D_4-D_5,\, -2D_5-D_6+D_7\}
\]
with respect to which the intersection matrix is 
\[
A_{\mathcal{B}'} = \begin{pmatrix} 0 & 1  \\ 1 & 0 & \multicolumn{2}{c}{\mbox{\Huge\boldmath$O$}}\\ & & -2 & 0 \\ \multicolumn{2}{c}{\mbox{\Huge\boldmath$O$}} & 0 & -4\end{pmatrix}. 
\]
Thus, $\Pic_{\!\Delta_1}=U\oplus\langle -2\rangle\oplus\langle -4\rangle$, which is clearly primitively embedded into the $K3$ lattice. 

\bigskip

\begin{figure}[!htb]
\begin{center}
\includegraphics[width=.4\linewidth]{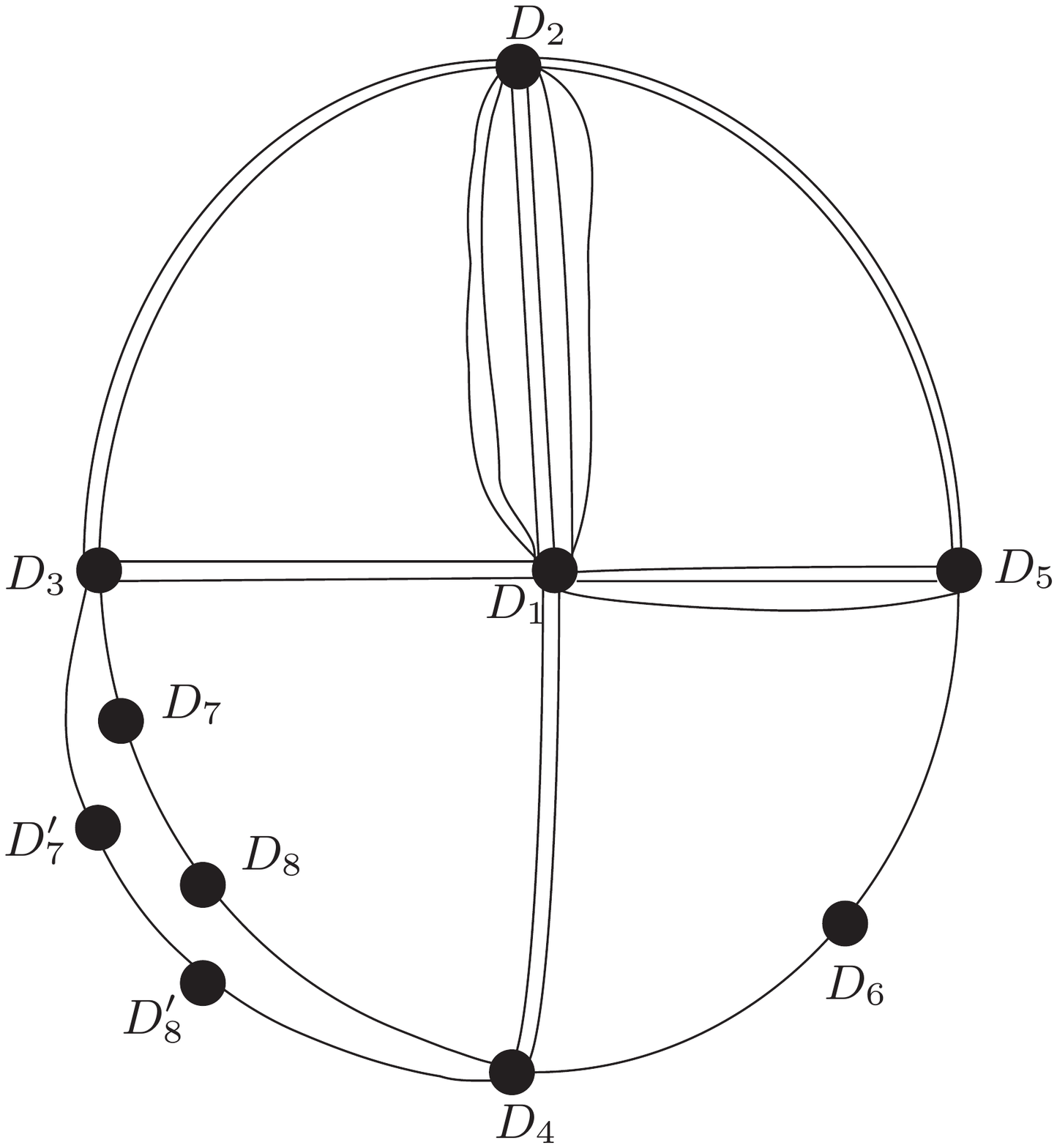}
\caption{Nodes are divisors on a general section in $\mathcal{F}_2$. 
There are $n$ lines connecting $D_i$ and $D_j$ if $D_i.D_j=n$. }
\label{Case2}
\end{center}
\end{figure}

We lebel the primitive vectors that generate rays of the fan defining the $3$-fold $\mathbb{P}_{\!\Delta_2}$, or equivalently, the lattice points in $\Delta_2^*$ as follows:
\[
\begin{array}{lll}
v_1=(0,0,1), & v_2=(-1,-1,-3), & v_3=(0,1,0),\\
v_4=(3,4,6), & v_5=(1,0,0), & v_6=(2,2,3), \\
v_7=(1,2,2), & v_8=(2,3,4). 
\end{array}
\]
Let $D_i:=\overline{{\rm Orb}(\mathbb{R}_{\geq 0}v_i)}|_{{-}K_{\mathbb{P}_{\!\Delta_2}}}$ for $i=1,\ldots ,6$ be restricted torus-invatiant divisors, and $D_7, \, D_7'$ and $D_8,\, D_8'$ be components of the divisor $\overline{{\rm Orb}(\mathbb{R}_{\geq 0}v_7)}|_{{-}K_{\mathbb{P}_{\!\Delta_2}}}$, and $\overline{{\rm Orb}(\mathbb{R}_{\geq 0}v_8)}|_{{-}K_{\mathbb{P}_{\!\Delta_2}}}$, respectively. 
One computes the self-intersection numbers by the formula (\ref{SelfIntersection}), 
\[
\begin{array}{lll}
D_1^2 = 2\cdot 8 -2 = 14, &
D_2^2 = 2\cdot 2 -2 = 2, &
D_3^2 = 2\cdot 0 -2 = -2, \\
D_4^2 = 2\cdot 0 -2 = -2, &
D_5^2 = 2\cdot 1 -2 = 0, &
D_6^2 = D_7^2 = D_7'^2 = D_8^2 = D_8'^2 = -2. 
\end{array}
\]
One also has the graph of intersections among these divisors by the formula (\ref{Intersection}) as in Figure \ref{Case2}. 

By solving the linear system
\[
\sum_{i=1}^8 ({}^tv_i, \mbi{e}_j) D_i= 0, 
\]
for $j=1,2,3$, where $\mbi{e}_j$ is the $j$-th column of the identity matrix of size $3$, and $(\mbi{x}, \mbi{y})$ is the inner product on $\mathbb{R}^3$, one obtains a set of linearly-independent divisors $\mathcal{B}=\{ D_4,\, D_5,\, D_6,\, D_7,\, D_7',\, D_8,\, D_8'\}$, the intersection matrix with respect to which is 
\[
A_{\mathcal{B}} = \begin{pmatrix}
 -2 & 0 & 1 & 0 & 0 & 1 & 1 \\
 0 & 0 & 1 & 0 & 0 & 0 & 0 \\
 1 & 1 & -2 & 0 & 0 & 0 & 0 \\
 0 & 0 & 0 & -2 & 0 & 1 & 0 \\
 0 & 0 & 0 & 0 & -2 & 0 & 1 \\
 1 & 0 & 0 & 1 & 0 & -2 & 0 \\
 1 & 0 & 0 & 0 & 1 & 0 & -2 \\ 
\end{pmatrix}. 
\]
By a translation $PA_{\mathcal{B}}{}^tP$ with 
\[
P = \begin{pmatrix}
 0 & 1 & 1 & 0 & 0 & 0 & 0 \\
 0 & 1 & 0 & 0 & 0 & 0 & 0 \\
 0 & 0 & 0 & 0 & 1 & 0 & 0 \\
 0 & 0 & 0 & 0 & 0 & 0 & 1 \\
 1 & -1 & 0 & 0 & 0 & 0 & 0 \\
 0 & 0 & 0 & 0 & 0 & 1 & 0 \\
 0 & 0 & 0 & 1 & 0 & 0 & 0 \\
\end{pmatrix}, 
\]
one gets a new basis 
\[
\mathcal{B}' = \{ D_5+D_6,\, D_5,\, D_7',\, D_8',\, D_4-D_5,\, D_8,\, D_7\}
\]
with respect to which the intersection matrix is 
\[
A_{\mathcal{B}'} = \begin{pmatrix}
 0 & 1  \\
 1 & 0 & \multicolumn{5}{c}{\mbox{\Huge\boldmath$O$}}\\
  &  & -2 & 1 & 0 & 0 & 0 \\
  &  & 1 & -2 & 1 & 0 & 0 \\
\multicolumn{2}{c}{\mbox{\Huge\boldmath$O$}}   & 0 & 1 & -2 & 1 & 0 \\
  &  & 0 & 0 & 1 & -2 & 1 \\
  &  & 0 & 0 & 0 & 1 & -2 \\
\end{pmatrix}. 
\]
Thus, $\Pic_{\!\Delta_2}=U\oplus A_5$, which is clearly primitively embedded into the $K3$ lattice. \\

\begin{figure}[!htb]
\begin{center}
\includegraphics[width=.4\linewidth]{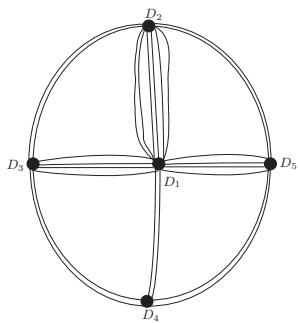}
\caption{Nodes are divisors on a general section in $\mathcal{F}_3$. 
There are $n$ lines connecting $D_i$ and $D_j$ if $D_i.D_j=n$.}
\label{Case3}
\end{center}
\end{figure}

(2)\, We lebel the primitive vectors that generate rays of the fan defining the $3$-fold $\mathbb{P}_{\!\Delta_3}$, or equivalently, the lattice points in $\Delta_3^*$ as follows:
\[
\begin{array}{lll}
v_1=(0,0,1), & v_2=(-1,-1,-3), & v_3=(0,1,0),\\
v_4=(1,1,1), & v_5=(1,0,0). 
\end{array}
\]
Let $D_i:=\overline{{\rm Orb}(\mathbb{R}_{\geq 0}v_i)}|_{{-}K_{\mathbb{P}_{\!\Delta_3}}}$ for $i=1,\ldots ,5$ be restricted torus-invatiant divisors. 
One computes the self-intersection numbers by the formula (\ref{SelfIntersection}), 
\[
\begin{array}{lll}
D_1^2 = 2\cdot 9 -2 = 16, &
D_2^2 = 2\cdot 3 -2 = 4, &
D_3^2 = 2\cdot 1 -2 = 0, \\
D_4^2 = 2\cdot 0 -2 = -2, &
D_5^2 = 2\cdot 1 -2 = 0. 
\end{array}
\]
One also has the graph of intersections among these divisors by the formula (\ref{Intersection}) as in Figure \ref{Case3}. 

By solving the linear system
\[
\sum_{i=1}^5 ({}^tv_i, \mbi{e}_j)D_i = 0, 
\]
for $j=1,2,3$, where $\mbi{e}_j$ is the $j$-th column of the identity matrix of size $3$, and $(\mbi{x}, \mbi{y})$ is the inner product on $\mathbb{R}^3$,  one obtains a set of linearly-independent divisors $\mathcal{B}=\{ D_4,\, D_5\}$, the intersection matrix with respect to which is 
\[
A_{\mathcal{B}} = \begin{pmatrix} -2 & 2\\ 2 & 0\end{pmatrix}. 
\]
By a translation $PA_{\mathcal{B}}{}^tP$ with 
\[
P = \begin{pmatrix} 1 & 0 \\ 1 & 1 \end{pmatrix}, 
\]
one gets a new basis 
\[
\mathcal{B}' = \{ D_4,\, D_4+D_5\}
\]
with respect to which the intersection matrix is 
\[
A_{\mathcal{B}'} = \begin{pmatrix} -2 & 0  \\ 0 & 2\end{pmatrix}. 
\]
Thus, $\Pic_{\!\Delta_3}=\langle -2\rangle\oplus\langle 2\rangle$, which is clearly primitively embedded into the $K3$ lattice. 

\begin{figure}[!htb]
\begin{center}
\includegraphics[width=.5\linewidth]{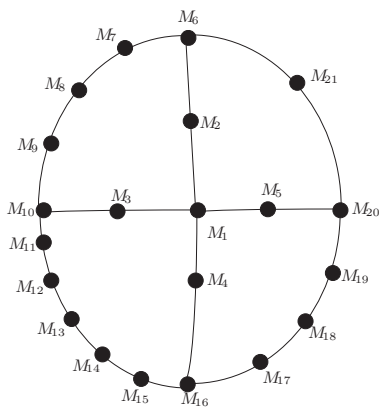}
\caption{Nodes are divisors on a general section in $\mathcal{F}_3'$. 
There are $n$ lines connecting $D_i$ and $D_j$ if $D_i.D_j=n$.}
\label{Case3Dual}
\end{center}
\end{figure}

(3)\, We lebel the primitive vectors that generate rays of the fan defining the $3$-fold $\mathbb{P}_{\!\Delta_3^*}$, or equivalently, the lattice points in $\Delta_3$ as follows:
\[
\begin{array}{lll}
m_1= (-1, -1, 1), & m_2=(0, -1, 0), & m_3=(2,-1,0), \\
m_4=(-1,2,0), & m_5=(-1,0,0), & m_6=(1,-1,-1), \\
m_7=(2,-1,-1), & m_8=(3,-1,-1), & m_9=(4,-1,-1), \\
m_{10}=(5,-1,-1), & m_{11}=(4,0,-1), & m_{12}=(3,1,-1), \\
m_{13}=(2,2,-1), & m_{14}=(1,3,-1), & m_{15}=(0,4,-1), \\
m_{16}=(-1,5,-1), & m_{17}=(-1,4,-1), & m_{18}=(-1,3,-1), \\
m_{19}=(-1,2,-1), & m_{20}=(-1,1,-1), & m_{21}=(0,0,-1). 
\end{array}
\]
Let $M_i:=\overline{{\rm Orb}(\mathbb{R}_{\geq 0}m_i)}|_{{-}K_{\mathbb{P}_{\!\Delta_3^*}}}$ for $i=1,\ldots ,21$ be restricted torus-invatiant divisors. 
One computes the self-intersection numbers by the formula (\ref{SelfIntersection}), 
\[
\begin{array}{lll}
M_1^2 = 2\cdot 1 -2 = 0, &
M_i^2 = 2\cdot 0 -2 = -2, & \forall i=2,\ldots ,21. 
\end{array}
\]
One also has the graph of intersection numbers among these divisors by the formula (\ref{Intersection}) as in Figure \ref{Case3Dual}. 

By solving the linear system
\[
\sum_{i=1}^{21} ({}^tm_i, \mbi{e}_j)D_i = 0, 
\]
for $j=1,2,3$, where $\mbi{e}_j$ is the $j$-th column of the identity matrix of size $3$, and $(\mbi{x}, \mbi{y})$ is the inner product on $\mathbb{R}^3$,  one obtains a set of linearly-independent divisors $\mathcal{B}=\{ M_1, \cdots , M_{18}\}$, the intersection matrix with respect to which is 
\[
B_{\mathcal{B}} = \begin{pmatrix}
 0 & 1 & 1 & 1 & 1 & 0 & 0 & 0 & 0 & 0 & 0 & 0 & 0 & 0 & 0 & 0 & 0 & 0 \\
 1 & -2 & 0 & 0 & 0 & 1 & 0 & 0 & 0 & 0 & 0 & 0 & 0 & 0 & 0 & 0 & 0 & 0 \\
 1 & 0 & -2 & 0 & 0 & 0 & 0 & 0 & 0 & 1 & 0 & 0 & 0 & 0 & 0 & 0 & 0 & 0 \\
 1 & 0 & 0 & -2 & 0 & 0 & 0 & 0 & 0 & 0 & 0 & 0 & 0 & 0 & 0 & 1 & 0 & 0 \\
 1 & 0 & 0 & 0 & -2 & 0 & 0 & 0 & 0 & 0 & 0 & 0 & 0 & 0 & 0 & 0 & 0 & 0 \\
 0 & 1 & 0 & 0 & 0 & -2 & 1 & 0 & 0 & 0 & 0 & 0 & 0 & 0 & 0 & 0 & 0 & 0 \\
 0 & 0 & 0 & 0 & 0 & 1 & -2 & 1 & 0 & 0 & 0 & 0 & 0 & 0 & 0 & 0 & 0 & 0 \\
 0 & 0 & 0 & 0 & 0 & 0 & 1 & -2 & 1 & 0 & 0 & 0 & 0 & 0 & 0 & 0 & 0 & 0 \\
 0 & 0 & 0 & 0 & 0 & 0 & 0 & 1 & -2 & 1 & 0 & 0 & 0 & 0 & 0 & 0 & 0 & 0 \\
 0 & 0 & 1 & 0 & 0 & 0 & 0 & 0 & 1 & -2 & 1 & 0 & 0 & 0 & 0 & 0 & 0 & 0 \\
 0 & 0 & 0 & 0 & 0 & 0 & 0 & 0 & 0 & 1 & -2 & 1 & 0 & 0 & 0 & 0 & 0 & 0 \\
 0 & 0 & 0 & 0 & 0 & 0 & 0 & 0 & 0 & 0 & 1 & -2 & 1 & 0 & 0 & 0 & 0 & 0 \\
 0 & 0 & 0 & 0 & 0 & 0 & 0 & 0 & 0 & 0 & 0 & 1 & -2 & 1 & 0 & 0 & 0 & 0 \\
 0 & 0 & 0 & 0 & 0 & 0 & 0 & 0 & 0 & 0 & 0 & 0 & 1 & -2 & 1 & 0 & 0 & 0 \\
 0 & 0 & 0 & 0 & 0 & 0 & 0 & 0 & 0 & 0 & 0 & 0 & 0 & 1 & -2 & 1 & 0 & 0 \\
 0 & 0 & 0 & 1 & 0 & 0 & 0 & 0 & 0 & 0 & 0 & 0 & 0 & 0 & 1 & -2 & 1 & 0 \\
 0 & 0 & 0 & 0 & 0 & 0 & 0 & 0 & 0 & 0 & 0 & 0 & 0 & 0 & 0 & 1 & -2 & 1 \\
 0 & 0 & 0 & 0 & 0 & 0 & 0 & 0 & 0 & 0 & 0 & 0 & 0 & 0 & 0 & 0 & 1 & -2 \end{pmatrix}. 
\]
By a direct computation, we have 
\[
\rk B_{\mathcal{B}} = 18, \quad \sign B_{\mathcal{B}} =(1,17),\quad \det{B_{\mathcal{B}}} = -4. 
\]
It also immediately follws by considering a map $F:\mathbb{Z}^{19}\to\mathbb{Q}^{15}$ by 
\[
\begin{array}{l}
\mathtt{F[x_1, \, x_2, \, x_3, \, x_4, \, x_5, \, x_6, \, x_7, \, x_8, \, x_9,  \, x_{10}, \, x_{11}, \, x_{12}, \, x_{13}, \, x_{14}, \, x_{15}, \, x_{16},\, x_{17}, \, x_{18} \, s] := } \\
\hspace{3mm}\mathtt{ s (x_1, \, x_2, \, x_3, \, x_4, \, x_5, \, x_6, \, x_7, \, x_8, \, x_9,  \, x_{10}, \, x_{11}, \, x_{12}, \, x_{13}, \, x_{14}\, x_{15}, \, x_{16},\, x_{17}, \,  x_{18}).B_{\mathcal{B}}^{-1}}
\end{array}
\]
that there exist two distinct elements in the discriminant group of order $2$ that generate the group since all the coefficients of $F[\mathbf{x},2]$ are integral. 
Thus, the lattice $\Pic_{\!\Delta_3^*}$ has the discriminant group $\mathbb{Z}\slash 2\mathbb{Z} \oplus \mathbb{Z}\slash 2\mathbb{Z}$ that is isometric to the discriminant group of $\Pic_{\!\Delta_3}=\langle -2\rangle\oplus\langle 2\rangle$. 
Therefore, by the Lemma below, $U\oplus \Pic_{\!\Delta_3^*}$ is isometric to $\left(\Pic_{\!\Delta_3}\right)_{\Lambda_{K3}}^\perp$. 
\begin{lem}[Corollary 1.6.2~\cite{Nikulin80}]\label{orthogonal}
Let lattices $S$ and $T$ be primitively embedded into the $K3$ lattice. 
Then $S$ and $T$ are orthogonal to each other in the $K3$ lattice if and only if $q_S\simeq -q_T$, where $q_S$ (resp. $q_T$) is the discriminant form of $S$ (resp. $T$). 
\end{lem}

The assertion is proved. 
\QED

\remark 
Since $\rho_{\!\Delta_i} + \rho_{\!\Delta_i^*}\not=20$, for $i=1,2$ by Proposition \ref{L0Rho}, the families $\mathcal{F}_1=\mathcal{F}_{\!\Delta_1}$ and $\mathcal{F}_{\!\Delta_1^*}$, and $\mathcal{F}_2=\mathcal{F}_{\!\Delta_2}$ and $\mathcal{F}_{\!\Delta_2^*}$ are not lattice dual in the sense of Dolgachev \cite{DolgachevMirror}.

On the other hand, by the statement of part (2) of Theorem \ref{MainThm}, the families $\mathcal{F}_3$ and $\mathcal{F}_3'$ are lattice dual. 

\remark
Since the Picard lattices $\Pic_{\!\Delta_1}$ and $\Pic_{\!\Delta_2}$ both contain the hyperbolic lattice $U$ of rank $2$ as a sublattice, general sections of families $\mathcal{F}_1$ and $\mathcal{F}_2$ have a structure of Jacobian elliptic surface.

Since the Picard lattice $\Pic_{\!\Delta_3}$ contains divisors $D_5$ and $D_4$ satisfying $D_5^2=-2, \, D_4^2=0$, and $D_5.D_4=2$, general sections of the family has a structure of elliptic surface, but, not Jacobian. 
Indeed, there exist smooth curves $s$ represented by the divisor $D_5$, and $f$ by $D_4$ such that $s$ is genus $1$, $f$ is genus $0$, and $s.f = 2$, thus, this is a $2$-section. 

\section{Families of double sextic $K3$ surfaces branching along curves of $(2,3)$-torus type}
In this section, we introduce and study the families $\mathcal{F}_1,\, \mathcal{F}_2,$ and $\mathcal{F}_3$ that contain $K3$ surfaces that are double covering of $\mathbb{P}^2$ branching along curves of $(2,3)$-torus type. 
Since every Gorenstein $K3$ surface is birationally equivalent to a unique $K3$ surface, we identify these two surfaces. 

\begin{prop}
The family $\mathcal{F}_1$ contains $K3$ surfaces that are double covering of $\mathbb{P}^2$ branching along curves of $(2,3)$-torus type with singularities 
\[
\begin{array}{l}
A_{17}, \quad
A_2+A_{14}, \quad
A_5 + A_{11}, \quad
E_6 + A_{11}, \quad
2A_8, \quad
2A_2+A_{11},  \\
A_2+A_5+A_8, \quad
3A_5, \quad
3A_2+A_8, \quad
2A_2+2A_5, \textnormal{ and }
4A_2+A_5. 
\end{array}
\]
\end{prop}
\proof 

$\mathbf {[A_{17}]}$ \, 
The curve is defined by $F(X,Y,Z):=F_2^3+F_3^2 =0$~\cite{Pho} with 
\[
F_2(X,Y,Z) = YZ-X^2, \qquad
F_3(X,Y,Z) = -X^2Z + YZ^2 + Y^3. 
\]
The polynomial $F$ contains monomials $X^6$, $Y^6$, $X^4 Z^2$, $Y^4 Z^2$, $Y^3 Z^3$, and $Y^2 Z^4$. 
Thus the Gorenstein $K3$ surface $D(X\slash \mathbb{P}^2)$ is a section in $\mathcal{F}_1$. 
Moreover, by a direct computation ({\it e.g.} with Singular), it is verified that the curve has one singularity of type $A_{17}$ at $(0:0:1)$. 

$\mathbf{[A_2+A_{14}]}$ \,
The curve is defined by $F(X,Y,Z,t_5):=F_2^3+F_3^2 =0$~\cite{Pho} with 
\[
F_2(X,Y,Z) =  YZ-X^2, \qquad
F_3(X,Y,Z,t_5) = Y^3 + t_5 X Y^2 -X^2 Z+Y Z^2. 
\]
The polynomial $F$ contains monomials $X^6$, $Y^6$, $X^2 Y^4$, $X^4 Z^2$, $X Y^5$, $Y^4 Z^2$, $Y^3 Z^3$, and $Y^2 Z^4$. 
Thus the Gorenstein $K3$ surface $D(X\slash \mathbb{P}^2)$ is a section in $\mathcal{F}_1$. 
Moreover, by a direct computation ({\it e.g.} with Singular), it is verified that there exists a curve with singularities of type $A_{14}$ at $(0:0:1)$, and of type $A_2$ at $(-1:1:1)$. 
For instance, setting $t_5=1$ would do. 

$\mathbf{[A_5 + A_{11}]}$\, 
The curve is defined by $F(X,Y,Z,t_1,t_2,s):=F_2^3+F_3^2 =0$~\cite{Pho} with 
\begin{eqnarray*}
F_2(X,Y,Z) & = & YZ-X^2, \\
F_3(X,Y,Z,t_1,t_2,s) & = & -t_2X^3 -(t_1+t_2-s-1)Y^3 +X^2Y-t_1X^2Z \\
& &  +2(t_1+t_2-s-1)XY^2 + t_1YZ^2-(t_1+t_2-s)Y^2Z. 
\end{eqnarray*}
The polynomial $F$ contains monomials $X^6$, $Y^6$, $X^5 Y$,  $X^4 Y^2$, $X^3 Y^3$, $X^2 Y^4$, $X Y^5$, $X^3 Y^3$, $X^2 Y^4$, $X Y^5$, $X^4 Z^2$, $X^5 Z$, $Y^5 Z$, $Y^4 Z^2$, $Y^3 Z^3$,  and $Y^2 Z^4$. 
Thus the Gorenstein $K3$ surface $D(X\slash \mathbb{P}^2)$ is a section in $\mathcal{F}_1$. 
Moreover, by a direct computation ({\it e.g.} with Singular), it is verified that there exists a curve with singularities of type $A_{11}$ at $(0:0:1)$, and of type $A_5$ at $(1:1:1)$. 
For instance, setting $t_1=t_2=1$ and $s=0$ would do. 

$\mathbf{[E_6 + A_{11}]}$\,
The curve is defined by $F(X,Y,Z,t_1,t_2):=F_2^3+F_3^2 =0$~\cite{Pho} with 
\begin{eqnarray*}
F_2(X,Y,Z) & = & YZ-X^2, \\
F_3(X,Y,Z,t_1,t_2) & = & -t_2X^3+(-t_1-t_2+1)Y^3 +X^2Y  \\
& & +2(t_1+t_2-1)XY^2 + t_1YZ^2+(t_1+t_2)Y^2Z -t_1X^2Z.    
\end{eqnarray*}
The polynomial $F$ contains monomials $X^6$, $Y^6$, $X^5 Y$, $X^5 Z$, $X^4 Y^2$, $X^4 Z^2$, $X^3 Y^3$, $X^2 Y^4$, $X Y^5$, $Y^5 Z$, $Y^4 Z^2$,  $Y^3 Z^3$, and $Y^2 Z^4$.  
Thus the Gorenstein $K3$ surface $D(X\slash \mathbb{P}^2)$ is a section in $\mathcal{F}_1$. 
Moreover, by a direct computation ({\it e.g.} with Singular), it is verified that there exists a curve with singularities of type  $A_{11}$ at $(0:0:1)$, and of type $E_6$  at $(1:1:1)$. 
For instance, setting $t_1=-1$ and $t_2=1$ would do. 

$\mathbf{[2A_8]}$\,
The curve is defined by $F(X,Y,Z,t_1, t_2):=F_2^3+F_3^2 =0$~\cite{Pho} with 
\begin{eqnarray*}
F_2(X,Y,Z) & = & YZ-X^2, \\
F_3(X,Y,Z,t_1,t_2) & = &-(t_2+1)X^3+t_2Y^3+3t_2X^2Y \\
& & -3t_2XY^2-t_1X^2Z  +t_1YZ^2. 
\end{eqnarray*}
The polynomial $F$ contains monomials $X^6$, $Y^6$, $X^5 Y$, $X^4 Y^2$, $X^4 Z^2$, $X^3 Y^3$, $X^2 Y^4$, $X Y^5$, $X^5 Z$, $Y^4 Z^2$, and $Y^3 Z^3$. 
Thus the Gorenstein $K3$ surface $D(X\slash \mathbb{P}^2)$ is a section in $\mathcal{F}_1$. 
Moreover, by a direct computation ({\it e.g.} with Singular), it is verified that there exists a curve with singularities of type $A_{8}$ at $(0:0:1)$ and $(1:1:1)$. 
For instance, setting $t_1=t_2=1$ would do. 

$\mathbf{[2A_2+A_{11}]}$\,
The curve is defined by $F(X,Y,Z,t_4, t_5):=F_2^3+F_3^2 =0$~\cite{Pho} with 
\[
F_2(X,Y,Z)  =  YZ-X^2, \quad
F_3(X,Y,Z,t_4,t_5) = Y^3 + t_4 X^2 Y -X^2 Z+t_5 X Y^2 + YZ^2. 
\]
The polynomial $F$ contains monomials $X^6$, $Y^6$, $X^4 Y^2$, $X^4 Z^2$, $X^3 Y^3$, $X^2 Y^4$, $X Y^5$, $Y^4 Z^2$, $Y^3 Z^3$, and $Y^2 Z^4$. 
Thus the Gorenstein $K3$ surface $D(X\slash \mathbb{P}^2)$ is a section in $\mathcal{F}_1$. 
Moreover, by a direct computation ({\it e.g.} with Singular), it is verified that there exists a curve with singularities type $A_{11}$ at $(0:0:1)$, and $A_2$ at $(-1:1:1)$ and $(-2:4:1)$. 
For instance, setting $t_4=-2$ and $t_5=1$ would do.

$\mathbf{[A_2+A_5+A_8]}$\, 
The curve is defined by $F(X,Y,Z,t_1, t_2,s):=F_2^3+F_3^2 =0$~\cite{Pho} with 
\begin{eqnarray*}
F_2(X,Y,Z)  & =  & YZ-X^2, \\
F_3(X,Y,Z,t_1,t_2,s) & = & \left(-t_2 - \frac{23}{27}\right) X^3  - \left(t_2 - \frac{4}{27}\right) Y^3  \\
& & + \left(t_1 + t_2 + \frac{23}{27} - s\right) X^2Y + \left(t_2 - \frac{4}{27}\right) XY^2 \\
& & -t_1 X^2Z + XYZ  + t_1 Y Z^2 + \left(-t_1 - 1 + s\right) Y^2 Z. 
\end{eqnarray*}
The polynomial $F$ contains monomials  $X^6$, $Y^6$, $X^5 Y$, $X^4 Y^2$, $X^3 Y^3$, $X^2 Y^4$, $X Y^5$, $X^5 Z$, $X^4 Z^2$, $Y^5 Z$, $Y^4 Z^2$, $Y^3 Z^3$, and $Y^2 Z^4$. 
Thus the Gorenstein $K3$ surface $D(X\slash \mathbb{P}^2)$ is a section in $\mathcal{F}_1$. 
Moreover, by a direct computation ({\it e.g.} with Singular), it is verified that there exists a curve with singularities of type $A_8$ at $(0:0:1)$, of type $A_2$ at $(-1:1:1)$, and of type $A_5$ at $(1:1:1)$.
For instance, setting $t_1=t_2=s=1$ would do.

$\mathbf{[3A_5]}$\, 
The curve is defined by $F(X,Y,Z,t_1, t_2,t_3):=F_2^3+F_3^2 =0$~\cite{Pho} with 
\begin{eqnarray*}
F_2(X,Y,Z)  & =  & YZ-X^2, \\
F_3(X,Y,Z,t_1,t_2,t_3) & = & -t_3X^3 + X^2Y +\left(-t_1-\frac{t_2}{2}+t_3-1\right)X^2Z + t_3XYZ \\
& & \hspace{5mm} +\left(-\frac{t_2}{2}-1+t_3\right)Y^3 + (t_2+1-2t_3)Y^2Z +t_1YZ^2. 
\end{eqnarray*}
The polynomial $F$ contains monomials $X^6$, $Y^6$, $X^5 Y$, $X^5 Z$, $X^4 Y^2$, $X^4 Z^2$, $X^3 Y^3$, $X^2 Y^4$, $Y^5 Z$, $Y^4 Z^2$, $Y^3 Z^3$, and $Y^2 Z^4$. 
Thus the Gorenstein $K3$ surface $D(X\slash \mathbb{P}^2)$ is a section in $\mathcal{F}_1$. 
Moreover, by a direct computation ({\it e.g.} with Singular), it is verified that there exists a curve with singularities type $A_8$ at $(0:0:1)$, of type $A_{5}$ at $(0:0:1), (1:1:1)$, and $(-1:1:1)$ .
For instance, setting $t_1=t_2=t_3=1$ would do.

$\mathbf{[3A_2+A_8]}$\, 
The curve is defined by $F(X,Y,Z,t_3, t_4,t_5):=F_2^3+F_3^2 =0$~\cite{Pho} with 
\begin{eqnarray*}
F_2(X,Y,Z)  & =  & YZ-X^2, \\
F_3(X,Y,Z,t_3,t_4,t_5) & = & t_3 X^3+Y^3 +t_4 X^2Y-X^2 Z+YZ^2+t_5 X Y^2. 
\end{eqnarray*}
The polynomial $F$ contains monomials $X^6$, $Y^6$, $X^5 Y$, $X^5 Z$, $X^4 Y^2$, $X^4 Z^2$, $X^3 Y^3$, $X^2 Y^4$, $X Y^5$, $Y^4 Z^2$, $Y^3 Z^3$, and $Y^2 Z^4$. 
Thus the Gorenstein $K3$ surface $D(X\slash \mathbb{P}^2)$ is a section in $\mathcal{F}_1$. 
Moreover, by a direct computation ({\it e.g.} with Singular), it is verified that there exists a curve with singularities of type $A_{8}$ at $(0:0:1)$, and of type $A_2$ at $(-1:1:1), (-2:4:1)$, and $(2:4:1)$. 
For instance, setting $t_3=t_4=-4$ and $t_5=1$ would do.

$\mathbf{[2A_2+2A_5]}$\, 
The curve is defined by $F(X,Y,Z,t_1, t_2):=F_2^3+F_3^2 =0$~\cite{Pho} with 
\begin{eqnarray*}
F_2(X,Y,Z)  & =  & YZ-X^2, \\
F_3(X,Y,Z,t_1,t_2) & = & 3X^3 + Y^3 -(t_1+2)X^2Z +(t_1-t_2-1)X^2Y \\
& &  -3XY^2 + t_1YZ^2 +(-t_1+2+t_2)Y^2Z. 
\end{eqnarray*}
The polynomial $F$ contains monomials $X^6$, $Y^6$, $X^5 Y$, $X^5 Z$, $X^4 Y^2$, $X^4 Z^2$, $X^3 Y^3$, $X^2 Y^4$, $X Y^5$, $Y^5 Z$, $Y^4 Z^2$, and $Y^2 Z^4$.  
Thus the Gorenstein $K3$ surface $D(X\slash \mathbb{P}^2)$ is a section in $\mathcal{F}_1$. 
Moreover, by a direct computation ({\it e.g.} with Singular), it is verified that there exists a curve with singularities type $A_{5}$ at $(0:0:1)$ and $(1:1:1)$, and of type $A_2$ at $(-1:1:1)$ and $(2:4:1)$. 
For instance, setting $t_1=t_2=1$ would do.

$\mathbf{[4A_2+A_5]}$\,
The curve is defined by $F(X,Y,Z,t_2, t_3, t_4,t_5):=F_2^3+F_3^2 =0$~\cite{Pho} with 
\begin{eqnarray*}
F_2(X,Y,Z)  & =  & YZ-X^2, \\
F_3(X,Y,Z,t_2, t_3,t_4,t_5) & = & t_3 X^3+Y^3 +t_4 X^2Y+\left(t_2-1\right) X^2 Z \\
& & +t_5 X Y^2+YZ^2. 
\end{eqnarray*}
The polynomial $F$ contains monomials $X^6$, $Y^6$, $X^5 Y$, $X^5 Z$, $X^4 Y^2$, $X^4 Z^2$, $X^3 Y^3$, $X^2 Y^4$, $X Y^5$, $Y^4 Z^2$, $Y^3 Z^3$,  and $Y^2 Z^4$.  
Thus the Gorenstein $K3$ surface $D(X\slash \mathbb{P}^2)$ is a section in $\mathcal{F}_1$. 
Moreover, by a direct computation ({\it e.g.} with Singular), it is verified that there exists a curve with singularities of type $A_{5}$ at $(0:0:1)$, and of type $A_2$ at $(-1:1:1), (1:1:1), (-2:4:1)$, and $(2:4:1)$. 
For instance, setting $(t_2, t_3, t_4, t_5) = (4, 0,-5, 0)$ would do.

Thus the assetion is verified. 
\QED

\begin{prop}
The family $\mathcal{F}_2$ contains $K3$ surfaces that are double covering of $\mathbb{P}^2$ branching along curves of $(2,3)$-torus type with singularities 
\[
\begin{array}{l}
2A_5+E_6, \quad
A_5+2E_6, \quad
3E_6,  \quad
2A_2+A_5+E_6, \\
2A_2+2E_6, \textnormal{ and } \quad
A_2+E_6+A_{8}. 
\end{array}
\]
\end{prop}
\proof 

$\mathbf{[2A_5+E_6]}$\, 
The curve is defined by $F(X,Y,Z,t_2, t_3):=F_2^3+F_3^2 =0$~\cite{Pho} with 
\begin{eqnarray*}
F_2(X,Y,Z) & = & YZ-X^2, \\
F_3(X,Y,Z,t_2,t_3) & = & -t_3X^3 + X^2Y -\left(1+\frac{t_2}{2}-t_3\right)X^2Z + t_3XYZ \\
& & -\left(1+\frac{t_2}{2}-t_3\right)Y^3 + (1+t_2-2t_3)Y^2Z. 
\end{eqnarray*}
The polynomial $F$ contains monomials $X^6$, $Y^6$, $X^5 Y$, $X^5 Z$, $X^4 Y^2$, $X^4 Z^2$, $X^3 Y^3$, $X^2 Y^4$, $Y^5 Z$, $Y^4 Z^2$, and $Y^3 Z^3$. 
Thus the Gorenstein $K3$ surface $D(X\slash \mathbb{P}^2)$ is a section in $\mathcal{F}_2$. 
Moreover, by a direct computation ({\it e.g.} with Singular), it is verified that there exists a curve with singularities of type $E_6$ at $(0:0:1)$, and of type $A_{5}$ at $(1:1:1)$ and $(-1:1:1)$.  
For instance, setting $t_1=t_2=1$ would do.

$\mathbf{[A_5+2E_6]}$\, 
The curve is defined by $F(X,Y,Z,t_3):=F_2^3+F_3^2 =0$~\cite{Pho} with 
\begin{eqnarray*}
F_2(X,Y,Z) & = & YZ-X^2, \\
F_3(X,Y,Z,t_3) & = & -t_3X^3 -(1-t_3)Y^3 + X^2Y -(1-t_3)X^2Z \\
& & \hspace{5mm} + (1-2t_3)Y^2Z + t_3XYZ. 
\end{eqnarray*}
The polynomial $F$ contains monomials $X^6$, $Y^6$, $X^5 Y$, $X^5 Z$, $X^4 Y^2$, $X^3 Y^3$, $X^4 Z^2$, $X^2 Y^4$, $Y^3 Z^3$, and $Y^5 Z$. 
Thus the Gorenstein $K3$ surface $D(X\slash \mathbb{P}^2)$ is a section in $\mathcal{F}_2$. 
Moreover, by a direct computation ({\it e.g.} with Singular), it is verified that there exists a curve with singularities of type $A_{5}$ at $(0:0:1), (1:1:1)$ and $(-1:1:1)$. 
For instance, setting $t_3=-1$ would do.

$\mathbf{[3E_6]}$\, 
The curve is defined by $F(X,Y,Z):=F_2^3+F_3^2 =0$~\cite{Pho} with 
\[
F_2(X,Y,Z)  =  YZ-X^2, \qquad
F_3(X,Y,Z) = X^2Y -X^2Z -Y^3 + Y^2Z. 
\]
The polynomial $F$ contains monomials $X^6$, $Y^6$, $X^4 Y^2$, $X^4 Z^2$, $X^2 Y^4$, $Y^5 Z$, $Y^4 Z^2$, and $Y^3 Z^3$. 
Thus the Gorenstein $K3$ surface $D(X\slash \mathbb{P}^2)$ is a section in $\mathcal{F}_2$. 
Moreover, by a direct computation ({\it e.g.} with Singular), it is verified that the curve has singularities of type $E_{6}$ at $(0:0:1), (1:1:1)$ and $(-1:1:1)$. 

$\mathbf{[2A_2+A_5+E_6]}$\, 
The curve is defined by $F(X,Y,Z,t_2):=F_2^3+F_3^2 =0$~\cite{Pho} with 
\begin{eqnarray*}
F_2(X,Y,Z) & = & YZ-X^2, \\
F_3(X,Y,Z,t_2) & = &3X^3 + Y^3  -(1+t_2)X^2Y -3XY^2 -2X^2Z +(2+t_2)Y^2Z. 
\end{eqnarray*}
The polynomial $F$ contains monomials $X^6$, $Y^6$, $X^5 Y$, $X^5 Z$, $X^4 Y^2$, $X^4 Z^2$, $X^3 Y^3$, $X^2 Y^4$, $X Y^5$, $Y^5 Z$, $Y^4 Z^2$, and $Y^3 Z^3$. 
Thus the Gorenstein $K3$ surface $D(X\slash \mathbb{P}^2)$ is a section in $\mathcal{F}_2$. 
Moreover, by a direct computation ({\it e.g.} with Singular), it is verified that there exists a curve with singularities of type $E_{6}$ at $(0:0:1)$, of type $A_2$ at $(-1:1:1)$ and $(2:4:1)$, and of type $A_5$ at $(1:1:1)$. 
For instance, setting $t_2=1$ would do.

$\mathbf{[2A_2+2E_6]}$\, 
The curve is defined by $F(X,Y,Z):=F_2^3+F_3^2 =0$~\cite{Pho} with 
\[
F_2(X,Y,Z)  =  YZ-X^2, \quad
F_3(X,Y,Z) = 3X^3+ Y^3 -X^2Y -3XY^2  -2X^2Z  + 2Y^2Z. 
\]
The polynomial $F$ contains monomials $X^6$,  $Y^6$, $X^5 Y$, $X^5 Z$, $X^4 Y^2$, $X^4 Z^2$, $X^3 Y^3$, $X^2 Y^4$, $X Y^5$, $Y^5 Z$, $Y^4 Z^2$, and $Y^3 Z^3$. 
Thus the Gorenstein $K3$ surface $D(X\slash \mathbb{P}^2)$ is a section in $\mathcal{F}_2$. 
Moreover, by a direct computation ({\it e.g.} with Singular), it is verified that the curve has singularities of type $E_{6}$ at $(0:0:1)$ and $(1:1:1)$, and of type $A_2$ at $(-1:1:1)$ and $(2:4:1)$. 

$\mathbf{[A_2+E_6+A_{8}]}$\,  
The curve is defined by $F(X,Y,Z,t_1, t_2):=F_2^3+F_3^2 =0$~\cite{Pho} with 
\begin{eqnarray*}
F_2(X,Y,Z)  & =  & YZ-X^2, \\
F_3(X,Y,Z,t_1, t_2) & = & -\left(t_2+\frac{23}{27}\right) X^3+ \left(\frac{4}{27}-t_2\right) Y^3 \\
& & +\left(t_1+t_2+\frac{23}{27}\right) X^2 Y +\left(t_2-\frac{4}{27}\right) X Y^2 -t_1X^2 Z \\
& & -\left(t_1+1\right) Y^2 Z+t_1 Y Z^2+X Y Z. 
\end{eqnarray*}
The polynomial $F$ contains monomials $X^6$, $Y^6$, $X^5 Y$, $X^4 Y^2$, $X^3 Y^3$, $X^2 Y^4$, $X Y^5$, $X^5 Z$, $X^4 Z^2$, $Y^5 Z$, $Y^4 Z^2$, and $Y^3 Z^3$. 
Thus the Gorenstein $K3$ surface $D(X\slash \mathbb{P}^2)$ is a section in $\mathcal{F}_2$. 
Moreover, by a direct computation ({\it e.g.} with Singular), it is verified that there exists a curve with singularities of type $A_8$ at $(0:0:1)$, of type $A_2$ at $(-1:1:1)$, and of type $E_6$ at $(1:1:1)$. 
For instance, setting $t_1=t_2=1$ would do.

Thus the assetion is verified. 
\QED
\begin{prop}
The family $\mathcal{F}_3$ contains $K3$ surfaces that are double covering of $\mathbb{P}^2$ branching along the curve of $(2,3)$-torus type with singularities 
\[
4A_2+E_6. 
\]
The general $K3$ surfaces that are double covering of $\mathbb{P}^2$ branching along  the curve of $(2,3)$-torus type with singularities $6A_2$ cannot be contained in any of the three families, but in the full family of the weighted projective space with weight $(1,1,1,3)$. 
\end{prop}
\proof 

$\mathbf{[4A_2+E_6]}$\, 
The curve is defined by $F(X,Y,Z):=F_2^3+F_3^2 =0$~\cite{Pho} with 
\[
F_2(X,Y,Z)  =  YZ-X^2, \quad
F_3(X,Y,Z)  = Y^3 -X Y^2 -4 X Z^2-5 Y^2 Z+4 Y Z^2 +5 X Y Z. 
\]
The polynomial $F$ contains monomials $X^6$, $Y^6$, $X^2 Y^4$, $X^2 Z^4$, $X Y^5$, $Y^5 Z$, $Y^4 Z^2$, $Y^3 Z^3$, and $Y^2 Z^4$. 
Thus the Gorenstein $K3$ surface $D(X\slash \mathbb{P}^2)$ is a section in $\mathcal{F}_3$. 
Moreover, by a direct computation ({\it e.g.} with Singular), it is verified that the curve has singularities of type $E_6$ at $(1:1:1)$, and $A_2$ at $(0:0:1), (-1:1:1), (-2:4:1)$ and $(2:4:1)$. 

$\mathbf{[6A_2]}$\, 
The curve is defined by $F(X,Y,Z):=F_2^3+F_3^2 =0$~\cite{Pho} with 
\[
F_2(X,Y,Z)  =  YZ-X^2, \qquad
F_3(X,Y,Z)  = X^3+Y^3+Z^3. 
\]
The polynomial $F$ contains monomials $X^6$, $Y^6$,  $Z^6$, $X^3 Y^3$, $X^3 Z^3$, and $Y^3 Z^3$. 
Thus the Gorenstein $K3$ surface $D(X\slash \mathbb{P}^2)$ is a section neither in  $\mathcal{F}_1$ nor $\mathcal{F}_2$ nor $\mathcal{F}_3$, but in the full family of the weighted projective space with weight $(1,1,1,3)$. 
Moreover, by a direct computation ({\it e.g.} with Mathematica), it is verified that by B\'ezout's theorem, the conic $C_2:=(F_2=0)$ and the cubic $C_3:=(F_3=0)$ transversally intersect at the six points 
\begin{eqnarray*}
(e^{2/3i\pi}: e^{4/3i\pi}:1), 
(e^{2/3i\pi}\omega: e^{4/3i\pi}\omega^2:1), 
(e^{2/3i\pi}\omega^2: e^{4/3i\pi}\omega:1), \\
(e^{-2/3i\pi}: e^{-4/3i\pi}:1), 
(e^{-2/3i\pi}\omega^2: e^{-4/3i\pi}\omega:1), 
(e^{-2/3i\pi}\omega: e^{-4/3i\pi}\omega^2:1), 
\end{eqnarray*} 
at which points the curve $(F=0)$ has six singularities of type $A_2$. 
Note that the curve is the plane sextic that are studied by Zariski~\cite{Zariski29}. 

Let $(F'=0)$ be any deformation that preserves the type of singularity of the curve $(F=0)$. 
In this case, the polynomial $F'$ is a projective transformation of $F$. 
Thus, $F'$ also contains monomials $X^6,\, Y^6$, and $Z^6$. 
Thus, the $K3$ surface obtained by the double covering of $\mathbb{P}^2$ branching along the curve $(F'=0)$ cannot be contained in any families $\mathcal{F}_i,\, i=1,2,3$. 

The assertion is verified. 
\QED

Makiko Mase 

{\it Postal address 1}: \\
\begin{tabular}{r}
Department of Mathematics and Information Sciences \\
Tokyo Metropolitan University, \\
1-1 Hachioji-shi, Minami-osawa, Tokyo, Japan, 192-0397. 
\end{tabular} \\

{\it Postal address 2}: \\
\begin{tabular}{r}
Lehrstuhl f$\ddot{\textnormal{u}}$r Mathematik \\
Universit$\ddot{\textnormal{a}}$t Mannheim, \\
B6, 26 \quad 68131 Mannheim, Germany.  
\end{tabular} \\

{\it email address}: mtmase@arion.ocn.ne.jp
\end{document}